\documentclass[aap,nameyear,MSNbibl,seceqn,dvips]{arximspdf}
\usepackage{graphicx}

\doi{10.1214/10-AAP707}
\volume{21}
\issue{2}
\pubyear{2011}
\firstpage{645}
\lastpage{668}

\makeatletter
\newtheorem{theo}{Theorem}[section]
\newtheorem{cor}[theo]{Corollary}
\newtheorem{lemma}[theo]{Lemma}

\renewcommand{\epsilon}{\varepsilon}
\newcommand{\notag}{\nonumber}
\renewcommand{\cite}{\citet}
\renewcommand{\citep}[1]{[\citet{#1}]}
\newcommand{\var}{\operatorname{Var}}
\newcommand{\eqref}[1]{(\ref{#1})}
\makeatother

\begin{document}
\begin{frontmatter}

\title{On the transition from heavy traffic to heavy tails for the
M/G/1 queue: The regularly varying case}
\runtitle{Transition from heavy traffic to heavy tails}

\begin{aug}
\author[A]{\fnms{Mariana} \snm{Olvera-Cravioto}\corref{}\ead
[label=e1]{molvera@ieor.columbia.edu}},
\author[B]{\fnms{Jose} \snm{Blanchet}\ead[label=e2]{jose.blanchet@columbia.edu}}
\and
\author[C]{\fnms{Peter} \snm{Glynn} \ead[label=e3]{glynn@stanford.edu}}
\runauthor{M. Olvera-Cravioto, J. Blanchet and P. Glynn}
\affiliation{Columbia University, Columbia University and Stanford University}
\address[A]{M. Olvera-Cravioto\\Department of Industrial Engineering \\
and Operations Research \\ Columbia University \\ New York, New
York 10027 \\USA\\ \printead{e1}} %adresu isvedimo komanda gale!
\address[B]{J. Blanchet\\Department of Industrial Engineering \\
and Operations Research \\ Columbia University \\ New York, New
York 10027 \\USA\\ \printead{e2}}
\address[C]{P. Glynn\\Department of Management Science\\
and Engineering \\ Stanford University \\ Stanford, California
94305 \\USA\\ \printead{e3}}
\end{aug}

% HISTORY:
\received{\smonth{1} \syear{2009}}
\revised{\smonth{9} \syear{2010}}

% ABSTRACT
%
\begin{abstract}
Two of the most popular approximations for the distribution of the
steady-state waiting time, $W_\infty$, of the M/G/1 queue are the
so-called heavy-traffic approximation and heavy-tailed asymptotic,
respectively. If the traffic intensity, $\rho$, is close to 1 and
the processing times have finite variance, the heavy-traffic
approximation states that the distribution of $W_\infty$ is
roughly exponential at scale $O((1-\rho)^{-1})$, while the heavy tailed
asymptotic describes power law decay in the tail of the distribution of
$W_\infty$ for a fixed traffic intensity. In this
paper, we assume a regularly varying processing time distribution and
obtain a sharp threshold in terms of the tail value, or equivalently in
terms of $(1-\rho)$, that describes
the point at which the tail behavior transitions from the heavy-traffic
regime to the heavy-tailed asymptotic. We also provide new
approximations that are either uniform in the traffic intensity, or
uniform on the positive axis, that avoid the need to use different
expressions on the two regions defined by the threshold.
\end{abstract}

% KEYWORDS
%

\end{frontmatter}
%
%s1 ###
\section{Introduction}

A substantial literature has been developed over the last forty years
that recognizes the simplifications that arise in the analysis of
queueing systems in the presence of ``heavy traffic.'' The earliest
such ``heavy traffic'' approximation was that obtained by \citeauthor{Ki61}
(\citeyear{Ki61,Ki62}) for the steady-state waiting time $W_\infty$ for the G/G/1
queue. In particular, let $W_n$ be the waiting time (exclusive of
service) of the $n$th customer for a first-in first-out (FIFO)
single-server queue (with an infinite capacity waiting room) fed by a
renewal arrival process [with i.i.d. inter-arrival times $(\chi_n\dvtx
n\geq1)$] and an independent stream of i.i.d. processing times
$(V_n\dvtx
n\geq0)$. If $\rho\triangleq EV_1/E\chi_1 < 1$, then $W_n \Rightarrow
W_\infty$ as $n \to\infty$, where $W_\infty$ can be approximated via
%
%e1.1 ###
\begin{equation} \label{eq:HTrafficApprox}
W_\infty\stackrel{\mathcal{D}}{\approx} \frac{\var\chi_1 + \var
V_1}{2(E \chi_1 - E V_1)}   \operatorname{Exp}(1)
\end{equation}
when $\rho$ is close to 1. Here, $ \operatorname{Exp}(1)$ is an exponential r.v. with
mean one and $\stackrel{\mathcal{D}}{\approx}$ denotes ``has
approximately the same distribution as.'' A precise statement of the
limit theorem supporting the heavy traffic approximation \eqref
{eq:HTrafficApprox}
is given by \eqref{eq:HTrafficConvergence} below. The term ``heavy
traffic'' arises as a consequence of the fact that long queues tend to
form when $E V_1$ and $E \chi_1$ are roughly balanced.

The modern approach to justifying \eqref{eq:HTrafficApprox} involves
first showing that $(W_n\dvtx\break n\geq~0)$ can be approximated in heavy
traffic by a one-dimensional reflecting Brownian motion (RBM)
[see, e.g., \citeauthor{IgWh70a} (\citeyear{IgWh70a,IgWh70b})] and then verifying that the
steady-state r.v. $W_\infty$ can be approximated by that of the RBM
[\citeauthor{Sz90} (\citeyear{Sz90,Sz99})]. Similar methods apply, in significant generality, to
multi-station queueing networks. For example, \citet{Re84} proves a
functional limit theorem that justifies approximating a single class
multi-station queueing network by multi-dimensional RBM. Recent work of
\citet{GaZe05} establishes the associated steady-state convergence.
\citet{HaWi87} analyze the multi-dimensional RBM and show that it has
exponential tails.

On the other hand, if the processing times are heavy-tailed (e.g.,
regularly varying), there is a significant literature that establishes,
for various models, that the associated queueing system possesses a
heavy-tailed steady-state. A representative result of this type states
that when $\rho< 1$ for the G/G/1 FIFO queue described above (with
regularly varying processing times), we have
%
%e1.2 ###
\begin{equation} \label{eq:HTailApprox}
P(W_\infty> x) \sim\frac{\lambda}{1-\rho} \int_x^\infty P(V_1 > y) \,dy
\end{equation}
as $x \to\infty$, where $\lambda\stackrel{\triangle}{=} 1/E \chi_1$
[see, e.g., \citet{EmVe82}]. Corresponding heavy-tailed steady-state
asymptotics also exist in the context of queueing networks [see, e.g.,
\citet{BaScSc99} and \citet{BaFo04}].

At first, it may seem contradictory that the heavy-traffic theory
typically predicts exponential tails for the steady-state distribution,
whereas regularly varying heavy-tailed asymptotics predict power-law
decay in the steady-state tail. Of course, the key is to note that the
two families of results involve different types of limits, one as $\rho
\to1$ (heavy traffic) and the other as $x \to\infty$ (heavy tails).
The interesting mathematical issue here is therefore to send $\rho$ to
1 and $x \to\infty$ simultaneously, and to determine the $x$-value (as
a function of $\rho$) at which the steady-state distribution begins to
``feel'' the presence of the heavy tails in the processing times. In
particular, this paper develops a very explicit description, in the
setting of the M/G/1 queue (in which the arrival process is assumed
Poisson), of where the transition from the exponential heavy-traffic
approximation \eqref{eq:HTrafficApprox} to the heavy-tailed
approximation \eqref{eq:HTailApprox} occurs. As a corollary to our main
results (Corollary \ref{C.Threshold2}) we find that when the processing
times are regularly varying, then the tail probability $P(W_\infty>
x)$ sharply transitions at
%
%e1.3 ###
\begin{equation} \label{eq:transitionx}
x^* \approx\frac{1}{1-\rho} \log \biggl(\frac{1}{1-\rho} \biggr) \frac
{EV_1^2}{2 E V_1} (\alpha-2)
\end{equation}
from the approximation \eqref{eq:HTrafficApprox} to the approximation
\eqref{eq:HTailApprox} (where $\alpha$ is the tail index of the
regularly varying $V_1$). Roughly speaking, to the left of $x^*$, \eqref
{eq:HTrafficApprox} is valid whereas to the right of $x^*$, \eqref
{eq:HTailApprox} is appropriate. A companion paper [Olvera-Cravioto and Glynn (\citeyear{OlGl09})] provides uniform approximations
for $P(W_{\infty} > x)$ in the general subexponential case, and shows how in the setting of Weibullian
tails one can identify an intermediate zone in which neither the heavy-traffic asymptotic nor
the heavy-tailed asymptotic hold.

This result ties together two significant queueing theory literatures,
namely heavy traffic theory and heavy-tailed approximations. As the
first such result describing the transition from the heavy traffic
regime to the heavy-tailed asymptotic, it suggests the possibility of
similar such results for more complex systems and networks.
Furthermore, one of our main results, Theorem \ref{T.SumsApprox},
provides an approximation for the tail probability $P(W_\infty> x)$
that is uniform across all values of $\rho$, and that in numerical
experiments seems to perform very well. This new uniform approximation,
which takes advantage of the Pollaczek--Khintchine formula for the
M/G/1 queue, provides a significant numerical improvement over the
existing heavy-traffic and heavy-tail approximations that are commonly
used to approximate the tail of the r.v. $W_\infty$.

%%% The Main Results
%s2 ###
\section{The main results}

Let $(W_n(\rho)\dvtx  n\geq0)$ be the waiting time sequence for an M/G/1
FIFO queue that is fed by a Poisson arrival process having arrival rate
$\lambda= \rho/EV_1$ and independent i.i.d. processing times $(V_n\dvtx
n\geq0)$. We assume throughout the remainder of this paper (unless
otherwise noted) that $V_1$ has a regularly varying distribution with
tail index $\alpha> 2$, so that
\[
P(V_1 > x) \sim x^{-\alpha} L(x)
\]
as $x \to\infty$, where $L(\cdot )$ is slowly varying
 [see
page 412 of \citet{Asm2003}].

If $\rho< 1$, $W_n(\rho) \Rightarrow W_\infty(\rho)$ as $n \to\infty
$, where the Pollaczek--Khintchine formula  [see, e.g., page 237
of \citet{Asm2003}] guarantees that
%
%e2.1 ###
\begin{equation} \label{eq:PollKhin}
P\bigl(W_\infty(\rho) > \cdot \bigr) = \sum_{n=0}^\infty(1-\rho) \rho^n P(S_n
> \cdot ).
\end{equation}
Here, $S_n = X_1 + \cdots+ X_n$ (with $S_0=0$), where the $X_j$'s are
i.i.d. with common density $g(\cdot) = P(V_1 > \cdot )/EV_1$. The
heavy-tail result \eqref{eq:HTailApprox} translates, in the M/G/1
setting, into the asymptotic
%
%e2.2 ###
\begin{equation} \label{eq:HTailAsymptotic}
P\bigl(W_\infty(\rho) > x\bigr) \sim\frac{\rho}{1-\rho} P(X_1 > x),
\end{equation}
as $x \to\infty$. It is straightforward  [see, e.g., page 404
of \citet{Asm2003}] to show that \eqref{eq:HTailAsymptotic} in turn implies that
%
%e2.3 ###
\begin{equation} \label{eq:HTailAsymptotic2}
P\bigl(W_\infty(\rho) > x\bigr) \sim\frac{\lambda}{1-\rho} \cdot\frac
{x^{1-\alpha}}{\alpha-1} L(x)
\end{equation}
as $x \to\infty$.

Turning next to the heavy traffic limit theorem for $W_\infty(\rho)$
[due to \cite{Ki61}], its precise statement (in our M/G/1 setting) is that
%
%e2.4 ###
\begin{equation} \label{eq:HTrafficConvergence}
(1-\rho) W_\infty(\rho) \Rightarrow\frac{EV_1^2}{2EV_1}    \operatorname{Exp}(1)
\end{equation}
as $\rho\nearrow1$, providing theoretical support for the approximation
%
%e2.5 ###
\begin{equation} \label{eq:HTrafficApprox1}
P\bigl(W_\infty(\rho) > x\bigr) \approx\exp \bigl(-2(1-\rho) EV_1 x/E V_1^2 \bigr)
\end{equation}
when $\rho$ is close to 1. To get a sense of the point $x^* = x^*(\rho
)$ at which the heavy traffic approximation \eqref{eq:HTrafficApprox1}
transitions into the heavy-tail approximation \eqref
{eq:HTailAsymptotic2}, note that the point $x^*$ at which the
exponential \eqref{eq:HTrafficApprox1} crosses the power law tail \eqref
{eq:HTailAsymptotic2} must satisfy
%
%e2.6 ###
\begin{equation} \label{eq:match}
2 x^* (1-\rho) \frac{EV_1}{EV_1^2} \approx\log(1-\rho) + (\alpha-1)
\log x^*.
\end{equation}
This implies that $x^* \approx\kappa(1-\rho)^{-1} \log((1-\rho
)^{-1})$, where $\kappa= (\alpha-2) EV_1^2/(2 EV_1)$.

To make the above heuristic rigorous we look more closely at the
Pollaczek--Khintchine formula. First we note that the heavy-tail
asymptotic \eqref{eq:HTailAsymptotic} can be obtained by simply
substituting $P(S_n > \cdot)$ by $n P(X_1 > \cdot)$, that is, by using
the so-called subexponential asymptotic for $P(S_n > x)$. Such
asymptotics are typically stated for fixed values of $n$, but can be
shown to hold for $n \to\infty$ provided $n$ grows slowly enough
compared to $x$  [see, e.g., \citet{Bor00}; \citet{Roz89}]. In other words, we
can obtain the heavy-tail asymptotic from the first terms of \eqref
{eq:PollKhin},
\[
\sum_{n=1}^{N(x)} (1-\rho)\rho^n P(S_n > x) \approx\sum_{n=1}^{N(x)}
(1-\rho) \rho^n n P(X_1 > x) \sim\frac{\rho}{1-\rho} P(X_1 > x)
\]
for some appropriately defined $N(x)$.
This raises the question of whether we can also obtain the
heavy-traffic asymptotic directly from \eqref{eq:PollKhin}, and the
answer is yes. For large $n$, say $n \geq x/EX_1$, $P(S_n > x) = O(1)$,
so by simply replacing $P(S_n > x)$ by one we obtain
\[
\sum_{n= [x/EX_1]}^\infty (1-\rho)\rho^n P(S_n > x) \approx\sum_{n=
[x/EX_1]}^\infty (1-\rho)\rho^n = \rho^{[x/EX_1]} .
\]
Since as $\rho\nearrow1$, $\rho^{[x/EX_1]} \sim e^{-(1-\rho)x/ EX_1}
= e^{-2(1-\rho)EV_1 x/EV_1^2}$, we can recover the heavy-traffic
asymptotic from the last terms of \eqref{eq:PollKhin}.

This reasoning leads us to the observation that the transition of
$P(W_\infty> x)$ occurs at the level of the partial sums $P(S_n > x)$.
For the regularly varying case, the transition from the subexponential
asymptotic $nP(X_1 > x- nEX_1)$ to the CLT approximation $1- \Phi (
(x-nEX_1)/\sqrt{\var(X_1)} )$ [or its stable law counterpart when
$\var(X_1) = \infty$] occurs smoothly, which allows us to approximate
the Pollaczek--Khintchine formula directly and obtain an expression
that does not require $\rho$ to be close to one. Theorem \ref
{T.SumsApprox} below describes this (uniform in $\rho$) approximation,
and Theorem \ref{T.NewMain1} gives an equivalent formulation in terms
of more familiar asymptotic expressions. As corollaries, we obtain the
result regarding the transition from heavy-traffic to heavy-tail of
$P(W_\infty> x)$, both in terms of $\rho$ as a function of $x$ and $x$
as a function of $\rho$.

We also point out that similar versions of our results should also hold
for the GI/GI/1 case. The added difficulty lies in the fact that
although $W_\infty(\rho)$ still has a representation of the form
\[
W_\infty(\rho) = \sum_{n=1}^\infty(1-\theta) \theta^n P( Y_1 + \cdots+
Y_n > x) ,
\]
where the $Y_i$'s i.i.d regularly varying random variables
[see \citet{Asm2003}, Chapter~X.9], the distribution of the $Y_i$'s and the geometric
parameter $\theta$ are not explicitly known. In particular, both of
them depend on $\rho$, so a uniform in $\rho$ version of Theorem \ref
{OrigBor} and an asymptotic expression for $\theta(\rho)$ are required.
Such uniform in $\rho$ results have been recently developed in \citet
{BlGlLa09}. Proof techniques very similar to those given here can then
be used to obtain the GI/GI/1 equivalents of our results.

\begin{theo} \label{T.SumsApprox}
Suppose $P(V_1 > x) \sim L(x) x^{-\alpha}$ with $\alpha> 2$ and let
$\mu= EX_1 = EV_1^2/(2EV_1)$. Define $\beta= (2 \wedge(\alpha
-1))^{-1}$, $M(x) = \lfloor(x - x^\beta)/\mu\rfloor$, and
\[
S(\rho, x) = \sum_{n=1}^{M(x)} (1-\rho) \rho^n n P\bigl(X_1 > x-(n-1) \mu\bigr) .
\]
Then,
\[
\sup_{0 < \rho< 1}  \biggl| \frac{P(W_\infty(\rho) > x)}{S(\rho, x) +
\rho^{x/\mu}} - 1  \biggr| \to0
\]
as $x \to\infty$. Alternatively,
\[
\sup_{x \geq0}  \biggl| \frac{P(W_\infty(\rho) > x)}{S(\rho, x) + \rho
^{x/\mu}} - 1  \biggr| \to0
\]
as $\rho\nearrow1$.
\end{theo}

\begin{theo} \label{T.NewMain1}
Suppose $P(V_1 > x) \sim L(x) x^{-\alpha}$ with $\alpha> 2$ and let
$\mu= EX_1 = EV_1^2/(2EV_1)$ and $\gamma(x,\rho) = 1 - \rho^{x/\mu} -
\rho^{x/\mu} (1-\rho)x/\mu$. Then,
\[
\sup_{0 < \rho< 1}  \biggl| \frac{P(W_\infty(\rho) > x)}{ ({\rho
}/({1-\rho})) \gamma(x,\rho) P(X_1 > x) + \rho^{x/\mu}} - 1  \biggr| \to0
\]
as $x \to\infty$. Alternatively,
\[
\sup_{x \geq0}  \biggl| \frac{P(W_\infty(\rho) > x)}{ ({\rho}/({1-\rho}))
\gamma(x,\rho) P(X_1 > x) + \rho^{x/\mu}} - 1  \biggr| \to0
\]
as $\rho\nearrow1$.
\end{theo}

From Theorem \ref{T.NewMain1} we can derive the following corollary
stating the different regions where either the heavy-traffic
approximation or the heavy-tail asymptotic govern the tail behavior of
the steady-state waiting time. Corollary \ref{C.Threshold2} describes
the shape of the distribution of $W_\infty(\rho)$ for a fixed value of
$\rho$. On the other hand, Corollary \ref{C.Threshold1} can be of
practical use in understanding the sensitivity of a system to the
traffic intensity, since for a fixed value of $x$ it tells us how
$P(W_\infty(\rho) > x)$ changes as $\rho$ gets closer to one.

\begin{cor} \label{C.Threshold2}
Suppose $P(V_1 > x) \sim L(x) x^{-\alpha}$ with $\alpha> 2$ and let
$\kappa= (\alpha-2)EV_1^2/(2EV_1)$. Suppose that $y = y(\rho)$ satisfies
\[
y(\rho) = c\kappa(1-\rho)^{-1} \log \bigl( (1-\rho)^{-1} \bigr)
\]
for $\rho< 1$.
\begin{enumerate}[(b)]
\item[(a)] If $0 < c < 1$, then
%
%e2.7 ###
\begin{equation} \label{eq:HTrafficResult}
\sup_{0\leq x \leq y}  \biggl| \frac{P(W_\infty(\rho) > x)}{\exp
(-2(1-\rho)EV_1 x/EV_1^2 )} - 1 \biggr| \to0
\end{equation}
as $\rho\nearrow1$. Relation \eqref{eq:HTrafficResult} continues to
hold when $c = 1$, provided that $L(x)/(\log x)^{\alpha-1} \to0$ as $x
\to\infty$.

\item[(b)] If $c > 1$, then
%
%e2.8 ###
\begin{equation} \label{eq:HTailResult}
\sup_{x \geq y}  \biggl| \frac{P(W_\infty(\rho) > x)}{ ({\rho}/({1-\rho}))
P(X_1 > x)} - 1 \biggr| \to0
\end{equation}
as $\rho\nearrow1$. Relation \eqref{eq:HTailResult} continues to hold
when $c = 1$, provided that $L(x)/(\log x)^{\alpha-1} \to\infty$ as $x
\to\infty$.
\end{enumerate}
\end{cor}

The corresponding version in terms of $\rho$ as a function of $x$ is
given below.

\begin{cor} \label{C.Threshold1}
Suppose $P(V_1 > x) \sim L(x) x^{-\alpha}$ with $\alpha> 2$ and let
$\kappa= (\alpha-2)EV_1^2/(2EV_1)$. Suppose that $\hat\rho= \hat\rho
(x)$ satisfies
\[
\hat\rho(x) = 1 - c\kappa(\log x)/x.
\]
\begin{enumerate}[(b)]
\item[(a)] If $0 < c < 1$, then
%
%e2.9 ###
\begin{equation} \label{eq:HTrafficResult2}
\sup_{\hat\rho\leq\rho< 1}  \biggl| \frac{P(W_\infty(\rho) > x)}{\exp
 (-2(1-\rho)EV_1 x/EV_1^2 )} - 1 \biggr| \to0
\end{equation}
as $x \to\infty$. Relation \eqref{eq:HTrafficResult2} continues to
hold when $c = 1$, provided that $L(x)/\log x \to0$ as $x \to\infty$.
\item[(b)] If $c > 1$, then
%
%e2.10 ###
\begin{equation} \label{eq:HTailResult2}
\sup_{0 < \rho\leq\hat\rho}  \biggl| \frac{P(W_\infty(\rho) > x)}{
({\rho}/({1-\rho})) P(X_1 > x)} - 1 \biggr| \to0
\end{equation}
as $x \to\infty$. Relation \eqref{eq:HTailResult2} continues to hold
when $c = 1$, provided that $L(x)/\log x \to\infty$ as $x \to\infty$.
\end{enumerate}
\end{cor}

Note that Theorems \ref{T.SumsApprox} and \ref{T.NewMain1} suggest
different approximations for $P(W_\infty(\rho) > x)$. We tested both
approximations and found that
\[
H(\rho,x) = S(\rho,x) + \rho^{x/\mu}
\]
is better than its asymptotic counterpart and performs very well for
most values of $x$ and $\rho$. In Section \ref{S.Numerical} we analyze
how this approximation compares to using the simpler heavy-traffic and
heavy-tail asymptotics in the regions where they are valid, and we give
a couple of numerical examples.

It is instructive to contrast the behavior obtained in the above
regularly varying setting with what occurs in the light-tailed setting.
Suppose, in particular, that $E \exp(\theta V_1)< \infty$ for some
$\theta> 0$, and define $\theta^*(\rho)$ as the root of $\rho E \exp
(\theta^*(\rho)V_1) = 1$.

\begin{theo} \label{T.LightTails}
Suppose that $E \exp(\theta V_1) < \infty$ for some $\theta>0$.
\begin{enumerate}[(b)]
\item[(a)] If $y = y(\rho) = o ((1-\rho)^{-2} )$, then
\[
\sup_{0\leq x \leq y}  \biggl| \frac{P(W_\infty(\rho) > x)}{\exp
(-2(1-\rho)EV_1 x/EV_1^2 )} - 1 \biggr| \to0
\]
as $\rho\nearrow1$.

\item[(b)] For $x \geq0$,
\[
P\bigl(W_\infty(\rho) > x(1-\rho)^{-2}\bigr) \sim\exp \biggl(-2x(1-\rho)^{-1} \frac
{EV_1}{EV_1^2} + x \frac{EV_1^3}{3EV_1^2} - \frac{x}{4} \frac
{EV_1^2}{EV_1}  \biggr)
\]
as $\rho\nearrow1$.

\item[(c)] As $\rho\nearrow1$,
\[
\sup_{x \geq0}  \biggl| \frac{P(W_\infty(\rho) > x)}{\exp (-\theta
^*(\rho)x )} - 1 \biggr| \to0.
\]
\end{enumerate}
\end{theo}

Note that in contrast to the heavy-tailed setting, the heavy traffic
approximation is now valid over a larger range, namely up to tail
values of order $o ((1-\rho)^{-2} )$. At tail values of order
$(1-\rho)^{-2}$, the third moment of $V_1$ enters the asymptotic for
$P(W_\infty(\rho)>x)$ [see also \citet{Wh95} and \citet{BlGl07}].
Finally, part (c) shows that the Cram\'{e}r--Lundberg tail
asymptotic  [see, e.g., pages 365--369 of \citet{Asm2003}] is globally
valid in heavy traffic, showing the clear superiority of the Cram\'
{e}r--Lundberg asymptotic over the heavy traffic approximation when
$\rho$ is close to 1. On the other hand, for regularly varying tails,
any global approximation to $P (W_\infty(\rho) > \cdot  )$
must utilize both the heavy traffic approximation and the appropriate
tail asymptotic.

We close this section with a brief discussion of how the theory
described in this paper extends to the more general setting of
geometric random sums. Specifically, consider the random variable
\[
Z(p) = \sum_{i=1}^{N(p)} Y_i,
\]
where $(Y_i\dvtx  i \geq1)$ is a sequence of nonnegative nonlattice
i.i.d. random variables independent of the geometric r.v. $N(p)$ having
mass function
\[
P\bigl(N(p)=k\bigr) = (1-p)p^{k-1}
\]
for $k \geq1$ [see \citet{Kal1997} for various applied settings in
which such geometric random sums arise]. We assume that $Y_1$ is
regularly varying with finite variance, so that there exists $\beta>
2$ and a slowly varying function $L(\cdot)$ for which
\[
P(Y_1 > x) \sim x^{-\beta} L(x)
\]
as $x \to\infty$. Put $\tau= (\beta-1) E Y_1$.

\begin{theo} \label{T.Geom}
Let $\mu= EY_1$ and $\gamma(x,p) = 1 - (1-p)^{x/\mu} - (1- p)^{x/\mu}
px/\mu$. Then,
\[
\sup_{0 < p < 1}  \biggl| \frac{P(Z(p) > x)}{(1-p) \gamma(x,p) P(Y_1 >
x)/p + (1-p)^{x/\mu} } - 1  \biggr| \to0
\]
as $x \to\infty$. Alternatively,
\[
\sup_{x > 0}  \biggl| \frac{P(Z(p) > x)}{(1-p) \gamma(x,p) P(Y_1 > x)/p +
(1-p)^{x/\mu} } - 1  \biggr| \to0
\]
as $p \downarrow0$.
\end{theo}

\begin{cor} \label{C.GeomThres}
Suppose that $y = y(p)$ satisfies
\[
y(p) = c\tau p^{-1}\log (1/p )
\]
\begin{enumerate}[(b)]
\item[(a)] If $0<c<1$, then
%
%e2.11 ###
\begin{equation} \label{eq:GeoSum1}
\sup_{0\leq x \leq y}  \biggl| \frac{P(Z(p) > x)}{\exp(-px/EY_1)} - 1
 \biggr| \to0
\end{equation}
as $p \downarrow0$. Relation \eqref{eq:GeoSum1} continues to hold when
$c =1$, provided that $L(x)/(\log x)^{\beta} \to0$ as $x \to\infty$.

\item[(b)] If $c>1$, then
%
%e2.12 ###
\begin{equation} \label{eq:GeoSum2}
\sup_{x \geq y}  \biggl| \frac{P(Z(p) > x)}{P(Y_1 > x)/p} - 1  \biggr| \to0
\end{equation}
as $p \downarrow0$. Relation \eqref{eq:GeoSum2} continues to hold when
$c =1$, provided that $L(x)/(\log x)^{\beta} \to\infty$ as $x \to
\infty$.
\end{enumerate}
\end{cor}

%%% Proofs
%s3 ###
\section{Proofs}

In this section, we prove Theorems \ref{T.SumsApprox}, \ref{T.NewMain1}
and Corollary \ref{C.Threshold2}; the proofs of Theorem \ref{T.Geom}
and Corollary \ref{C.GeomThres} are essentially identical to those of
Theorem \ref{T.NewMain1} and Corollary \ref{C.Threshold2}. The proof of
Corollary \ref{C.Threshold1} is very similar in spirit to that of
Corollary \ref{C.Threshold2}, the difference being that it follows from
the uniform in $0 < \rho< 1$ statement of Theorem \ref{T.NewMain1}
instead of the uniform in $x > 0$. Theorem \ref{T.LightTails} follows
directly from Theorem 2 in \citet{BlGl07}.

We now turn our attention to the proof of Theorem \ref{T.SumsApprox}.
Recall that
%
%e3.1 ###
\begin{equation} \label{eq:Poll-Khin}
P\bigl(W_\infty(\rho) > x\bigr) = \sum_{n=0}^\infty(1-\rho) \rho^n P(S_n > x),
\end{equation}
where the $X_i$'s are i.i.d. with common density $g(\cdot ) = P(V_1 >
\cdot )/EV_1$ and $S_n = X_1 + \cdots+ X_n$.

Our analysis is based on the principle that we can approximate $P(S_n >
x)$ by the heavy tail asymptotic $n P(X_1 > x - (n-1) E[X])$ uniformly
in $n$ throughout the region of large deviations of $S_n$. Early
results of this kind are due to \citet{Nag82}, \citet{Roz89}, \cite{MikNag98},
\citet{Bor00}, and more recently, \citet
{DenDiekShn08}. The statement we present below is taken from \citet
{BorovBorov2008}, Theorems 3.4.1 and 4.4.1.

\begin{theo}[(Borovkov)] \label{OrigBor}
Let $Y_1, Y_2, \ldots$ be i.i.d. random variables having $E Y = 0$,
$\overline{F}(t) = P(Y_1 > t)$ and $\overline{F}(t) = t^{-\beta} L(t)$
where $L(\cdot)$ is slowly varying. Set $S_n = Y_1 + \cdots+ Y_n$, $n
\geq1$.
\begin{enumerate}[(b)]
\item[(a)] If $\beta> 2$ and $E Y^2 < \infty$, define $\sigma(n) = \sqrt
{(\beta-2)n\log n}$.
\item[(b)] If $\beta\in(1, 2)$ and $F(-t) \leq c \overline{F}(t)$ for $t >
0$ and some constant $c > 0$, define $\sigma(n) = \overline{F}^{-1}(1/n)$.
\end{enumerate}
Then, there exists a function $\varphi(t) \downarrow0$ as $t \uparrow
\infty$ such that
\[
\sup_{y \geq t \sigma(n)}  \biggl| \frac{P(S_n > y)}{nP(Y_1 > y)} -
1 \biggr| \leq\varphi(t)
\]
uniformly in $n$.
\end{theo}

Below we give an application of Borovkov's result to our particular setting.

\begin{lemma} \label{L.NewBorov}
Let $X_1, X_2, \ldots$ be i.i.d. nonnegative random variables with $\mu
= E[X] < \infty$, and $P(X_1 > t) = t^{-\alpha+1} L(t)$ where $L(\cdot
)$ is slowly varying and $\alpha> 2$. Set $S_n = X_1 + \cdots+ X_n$,
$n\geq1$. For any $(2 \wedge(\alpha-1))^{-1} < \gamma< 1$ define
$M_\gamma(x) = \lfloor(x - x^{\gamma})/\mu\rfloor$. Then, there
exists a function $\varphi(t) \downarrow0$ as $t \uparrow\infty$ such that
\[
\sup_{1 \leq n \leq M_\gamma(x)}  \biggl| \frac{P(S_n > x)}{n P(X_1 > x -
(n-1)\mu)} - 1  \biggr| \leq\varphi(x).
\]
\end{lemma}

\begin{pf}
Suppose first that $\alpha> 3$ and let $\sigma(n) = \sqrt{(\alpha
-2)n\log n}$. Since
\[
\frac{P(S_n > x)}{n P(X_1 > x-(n-1)\mu)} = \frac{P(S_n^* > x - n\mu)}{n
P(Y_1 > x - n\mu)} ,
\]
where $Y_i = X_i - \mu$ and $S_n^* = Y_1 + \cdots+ Y_n$. Then the
result will follow from Theorem \ref{OrigBor}(a) once we show that $(x
- n\mu)/ \sigma(n) \to\infty$ uniformly for $1 \leq n \leq M_\gamma
(x)$. To see this simply note that
\[
\frac{x-n\mu}{\sigma(n)} \geq\frac{x - M_\gamma(x) \mu}{\sigma
(M_\gamma(x))} \sim\sqrt{\frac{\mu}{\alpha-2}} \cdot\frac{x^{\gamma-
1/2}}{\sqrt{\log x }} .
\]
Since $\gamma> 1/2$, the above converges to infinity.

Suppose now that $\alpha\in(2,3)$ and note that $P(Y_1 \leq-t) = 0$
for $t \geq\mu$. Note also that since $\overline{F}(t) = P(Y_1 > t)$
is regularly varying with index $\alpha-1$, then $\sigma(n) = \overline
{F}^{-1}(1/n) = n^{1/(\alpha-1)} \tilde L(n)$ for some slowly varying
function $\tilde L(\cdot)$  [see \citet{BiGoTe1987}]. Then the result
will follow from Theorem \ref{OrigBor}(b) once we show that $(x - n\mu
)/ \sigma(n) \to\infty$ uniformly for $1 \leq n \leq M_\gamma(x)$. To
see this note that
\[
\frac{x-n\mu}{\sigma(n)} \geq\frac{x - M_\gamma(x) \mu}{\sigma
(M_\gamma(x))} \sim\frac{x^\gamma}{\sigma(x/\mu)} \sim\frac{x^{\gamma
- 1/(\alpha-1)}}{\mu^{-1/(\alpha-1)} \tilde L( x)} ,
\]
and since $\gamma> 1/(\alpha-1)$ the above converges to infinity.

The case $\alpha= 3$ is rather technical and does not provide
additional insights. We refer the reader to the internet supplement
\cite{OlBlGlSuppl09} for the details.
\end{pf}

We now give a lemma that will allow us to transform the statements of
the main results from being uniform in $0 < \rho< 1$ to being uniform
in $x > 0$, under the limiting regimes $x \to\infty$ and $\rho
\nearrow1$, respectively.

\begin{lemma} \label{L.Equivalence}
Suppose that
\[
\sup_{0 < \rho< 1}  \biggl| \frac{P(W_\infty(\rho) > x)}{A(\rho,x)} - 1
 \biggr| \to0
\]
as $x \to\infty$, where $A(\rho,x)$ satisfies
\[
\sup_{0 < x < (1-\rho)^{-\eta}}  | A(\rho,x) - 1  | \to0
\]
as $\rho\nearrow1$ for some $0 < \eta< 1$. Then,
\[
\sup_{x > 0}  \biggl| \frac{P(W_\infty(\rho) > x)}{A(\rho,x)} - 1  \biggr|
\to0
\]
as $\rho\nearrow1$.
\end{lemma}

\begin{pf}
We argue by contradiction. Suppose that there exists an $\epsilon> 0$
and a function $x\dvtx  (0,1) \to(0,\infty)$ such that $x(\phi) \geq(1-\phi
)^{-\eta}$ and
\[
 \biggl| \frac{P(W_\infty(\phi) > x(\phi))}{A(\phi,x(\phi))} -1  \biggr| >
\epsilon
\]
for all $0 < \phi< 1$. Then,
\[
\sup_{0 < \rho< 1}  \biggl| \frac{P(W_\infty(\rho) > x(\phi))}{A(\rho
,x(\phi))} -1  \biggr| \geq  \biggl| \frac{P(W_\infty(\phi) > x(\phi
))}{A(\phi,x(\phi))} -1  \biggr| > \epsilon.
\]
But this cannot be since by assumption,
\[
\lim_{\phi\nearrow1} \sup_{0 < \rho< 1}  \biggl| \frac{P(W_\infty(\rho
) > x(\phi))}{A(\rho,x(\phi))} -1  \biggr| = 0.
\]
It follows that
\[
\sup_{x \geq(1-\rho)^{-\eta}}  \biggl| \frac{P(W_\infty(\rho) >
x)}{A(\rho,x)} -1  \biggr| \to0
\]
as $\rho\nearrow1$. For $0 < x < (1-\rho)^{-\eta}$ note that
\begin{eqnarray*}
&&\lim_{\rho\nearrow1} \sup_{0 < x < (1-\rho)^{-\eta}}  \biggl| \frac
{P(W_\infty(\rho) > x)}{A(\rho,x)} -1  \biggr| \\
&& \qquad \leq\lim_{\rho\nearrow1} \sup_{0 < x < (1-\rho)^{-\eta}} \frac{
| P(W_\infty(\rho) > x) - 1  |}{A(\rho,x)} + \lim_{\rho\nearrow
1} \sup_{0 < x < (1-\rho)^{-\eta}} \frac{ | A(\rho,x) - 1
|}{A(\rho,x)} \\
&& \qquad = \lim_{\rho\nearrow1} \sup_{0 < x < (1-\rho)^{-\eta}}  \bigl|
P\bigl(W_\infty(\rho) > x\bigr) - 1  \bigr| .
\end{eqnarray*}
The last limit is zero by the standard heavy traffic limit.
\end{pf}

Throughout the rest of this section let $\mu= EX_1$, $\beta= (2
\wedge(\alpha-1))^{-1}$, $M(x) = \lfloor(x-x^\beta)/\mu\rfloor$, and
\[
S(\rho, x) = \sum_{n=1}^{M(x)} (1-\rho) \rho^n n P\bigl(X_1 > x-(n-1) \mu\bigr) .
\]

\begin{pf*}{Proof of Theorem \ref{T.SumsApprox}}
We will prove the uniform in $0 < \rho< 1$ asymptotic, since the
statement regarding the uniformity in $x > 0$ will follow from Lemma
\ref{L.Equivalence} by noting that
\begin{eqnarray*}
\sup_{0 < x < (1-\rho)^{-1/4}} S(\rho,x) &\leq&\sup_{0 < x < (1-\rho
)^{-1/4}} \frac{(1-\rho) M(x) (M(x)+1)}{2}\\
& \leq&\frac{(1-\rho
)^{1/2}}{\mu^2},
\end{eqnarray*}
which clearly converges to zero. Throughout the proof $C > 0$ is a
generic constant.

Fix $\beta< \gamma< 1 \wedge\beta(\alpha-1)$ and define $M_\gamma(x)
= \lfloor(x - x^\gamma)/\mu\rfloor$. Then, by Lemma~\ref{L.NewBorov},
there exists a function $\varphi_1(t) \to0$ as $t \to\infty$ such that
\[
\sup_{1\leq n \leq M_\gamma(x)}  \biggl| \frac{P(S_n > x)}{n P (X_1 >
x - (n-1)\mu )} - 1  \biggr| \leq\varphi_1(x).
\]
By \eqref{eq:Poll-Khin} we have
\begin{eqnarray*}
&&\bigl| P\bigl(W_\infty(\rho) > x\bigr) - S(\rho, x) - \rho^{x/\mu}  \bigr| \\
&& \qquad \leq\varphi_1(x) \sum_{n=1}^{M_\gamma(x)} (1-\rho) \rho^n nP\bigl(X_1 >
x-(n-1)\mu\bigr) \\
&& \qquad  \quad {}  +   \sum_{n=M_\gamma(x)+1}^{M(x)} (1-\rho)
\rho^n nP\bigl(X_1 > x - (n-1)\mu\bigr)\\
&& \qquad  \quad {}  +  \sum_{n= M_\gamma(x)+1}^{\lfloor x/\mu
\rfloor} (1-\rho) \rho^n P(S_n > x)\\
&& \qquad  \quad {} +  \Biggl| \sum_{n = \lfloor x/\mu
\rfloor+ 1}^\infty (1-\rho) \rho^n P(S_n > x) - \rho^{x/\mu}  \Biggr|.
\end{eqnarray*}
Clearly,
%
%e3.2 ###
\begin{equation} \label{eq:Term1}
\varphi_1(x) \sum_{n=1}^{M_\gamma(x)} (1-\rho) \rho^n nP\bigl(X_1 >
x-(n-1)\mu\bigr) \leq\varphi_1(x) S(\rho,x).
\end{equation}
Fix $0 < \epsilon< \min\{\alpha- 2, (\beta(\alpha-1)-\gamma)/\beta\}
$. The second term is bounded by
\begin{eqnarray*}
&& \sum_{n=M_\gamma(x)+1}^{M(x)} (1-\rho) \rho^n nP\bigl(X_1 > x - (n-1)\mu\bigr)
\\
&& \qquad \leq C x (1-\rho) \sum_{n=M_\gamma(x)+1}^{M(x)} \rho^n \bigl(x - (n-1)\mu
\bigr)^{-\alpha+1+\epsilon} .
\end{eqnarray*}
Since $g(n) = \rho^n (x-(n-1)\mu)^{-\alpha+1+\epsilon}$ is convex in $n$,
%
%e3.3 ###
\begin{eqnarray}\label{eq:ConvexBd}
&& x (1-\rho) \sum_{n= M_\gamma(x)+1}^{M(x)} g(n) \notag\\
&& \qquad \leq x (1-\rho) \bigl(M(x) - M_\gamma(x)\bigr) \max\bigl\{ g\bigl(M_\gamma(x)+1\bigr), g(M(x))
\bigr\} \nonumber
\\[-8pt]
\\[-8pt]
&& \qquad \leq C (1-\rho) x^{1+\gamma} \max \bigl\{ \rho^{(x-x^\gamma)/\mu}
x^{-\gamma(\alpha-1-\epsilon)},
\rho^{(x-x^\beta)/\mu- 1} x^{-\beta(\alpha-1-\epsilon)}  \bigr\}
\nonumber\\
&& \qquad \leq C (1-\rho) x^{1+\gamma-\beta(\alpha-1-\epsilon)} \rho^{(x-x^\gamma
)/\mu},\notag
\end{eqnarray}
where our choice of $\epsilon$ guarantees that $\gamma-\beta(\alpha
-1-\epsilon) < 0$. Also, we have
%
%e3.4 ###
\begin{eqnarray} \label{eq:BoundByOne}
\sum_{n= M_\gamma(x)+1}^{\lfloor x/\mu\rfloor} (1-\rho) \rho^n P(S_n >
x) &\leq&\rho^{M_\gamma(x)+1}  \bigl( 1 - \rho^{\lfloor x/\mu\rfloor-
M_\gamma(x)}  \bigr) \nonumber
\\[-8pt]
\\[-8pt]
&\leq& C \rho^{(x-x^\gamma)/\mu} x^\gamma|\log\rho|.\nonumber
\end{eqnarray}
To derive the last bound let $K(x) = \lfloor(x+ x^\gamma)/\mu\rfloor
$. Then,
\begin{eqnarray*}
&& \Biggl| \sum_{n = \lfloor x/\mu\rfloor+ 1}^\infty (1-\rho) \rho^n
P(S_n > x) - \rho^{x/\mu}  \Biggr| \\
&& \qquad = \rho^{x/\mu} - \rho^{\lfloor x/\mu\rfloor+1} + \sum_{n=\lfloor
x/\mu\rfloor+ 1}^\infty(1-\rho) \rho^n P(S_n \leq x) \\
&& \qquad \leq\rho^{x/\mu}  (1 - \rho ) + \sum_{n=\lfloor x/\mu
\rfloor+ 1}^{K(x)} (1-\rho)\rho^n + \sum_{n=K(x) + 1}^\infty(1-\rho
)\rho^n P( S_n \leq x) .
\end{eqnarray*}
It is easy to check that
\[
\sum_{n=\lfloor x/\mu\rfloor+ 1}^{K(x)} (1-\rho)\rho^n \leq C \rho
^{x/\mu} x^\gamma|\log\rho| .
\]
For the tail of the sum let $Y_i = \mu- X_i$ and $S_n^* = Y_1 + \cdots
+ Y_n$. Let $b_n$ be the scaling for which $Z_n = S_n^*/b_n \Rightarrow
Z$, where $Z$ is a stable random variable. Note that $b_n = n^{\beta}
L_0(n)$ for some slowly varying $L_0(\cdot)$. It follows that for all
$n > K(x)$,
\[
P(S_n \leq x) = P \biggl(Z_n \geq\frac{n\mu- x}{b_n}  \biggr) \leq P
\biggl(Z_n \geq\frac{(K(x)+1)\mu- x}{b_{K(x)+1}}  \biggr) ,
\]
where
\[
\frac{(K(x)+1)\mu- x}{b_{K(x)+1}} \geq\frac{\mu^\beta x^\gamma}{ (x +
x^\gamma)^\beta L_0((x-x^\gamma)/\mu)} \geq\frac{c x^{\gamma-\beta}}{L_0(x)}
\]
for some constant $c > 0$. It follows that
\[
\sup_{n > K(x)} P(S_n \leq x) \leq\sup_{n > K(x)} P\bigl(Z_n \geq c
x^{\gamma-\beta}/L_0(x)\bigr) \leq\varphi_2(x)
\]
for some $\varphi_2(t) \to0$ as $t \to\infty$. Hence,
\[
\sum_{n=K(x) + 1}^\infty(1-\rho)\rho^n P( S_n \leq x) \leq\varphi
_2(x) \sum_{n=K(x) + 1}^\infty(1-\rho)\rho^n \leq\varphi_2(x) \rho
^{x/\mu} .
\]
We thus have that
%
%e3.5 ###
\begin{eqnarray}\label{eq:Term4}
&& \Biggl| \sum_{n = \lfloor x/\mu\rfloor+ 1}^\infty (1-\rho) \rho^n
P(S_n > x) - \rho^{x/\mu}  \Biggr| \nonumber
\\[-8pt]
\\[-8pt]
&& \qquad \leq\varphi_2(x) \rho^{x/\mu} + C \rho^{x/\mu} |\log\rho| x^\gamma.\nonumber
\end{eqnarray}
Combining \eqref{eq:Term1}--\eqref{eq:Term4} gives
\begin{eqnarray*}
& &\bigl| P\bigl(W_\infty(\rho) > x\bigr) - S(\rho,x) - \rho^{x/\mu}  \bigr| \\
&& \qquad \leq\varphi_1(x) S(\rho,x) + \varphi_2(x) \rho^{x/\mu} + C \rho
^{(x-x^\gamma)/\mu} |\log\rho| x^\upsilon,
\end{eqnarray*}
where $\upsilon= \max\{1+\gamma-\beta(\alpha-1-\epsilon), \gamma\} \in
(0,1)$. It only remains to show that $\rho^{(x-x^\gamma)/\mu} |\log\rho
| x^\upsilon= o  (S(\rho,x)+\rho^{x/\mu} )$ uniformly in $0 <
\rho< 1$. To see this let $\rho(x) = 1 - (x^\upsilon\log x)^{-1}$, then
\begin{eqnarray*}
\sup_{\rho(x) \leq\rho< 1} \frac{\rho^{(x-x^\gamma)/\mu} |\log\rho|
x^\upsilon}{S(\rho,x) + \rho^{x/\mu}} &\leq&\sup_{\rho(x) \leq\rho<
1} e^{|\log\rho| x^\gamma/ \mu} |\log\rho| x^\upsilon \\
&=& e^{|\log\rho(x)| x^\gamma/ \mu} |\log\rho(x)| x^\upsilon\\
&\leq&\frac{C}{\log x} \to0 ,
\end{eqnarray*}
and since $S(\rho,x) \geq P(X_1 > x) \sum_{n=1}^{\lfloor x/\mu\rfloor}
(1-\rho) \rho^n$,
\begin{eqnarray*}
&&\sup_{0 < \rho< \rho(x) } \frac{\rho^{(x-x^\gamma)/\mu} |\log\rho|
x^\upsilon}{S(\rho,x) + \rho^{x/\mu}} \\
&& \qquad \leq\sup_{0 < \rho< \rho(x) } \frac{\rho^{(x-x^\gamma)/\mu} |\log\rho
| x^{\upsilon}}{L(x) x^{-\alpha+1} \rho(1- \rho^{\lfloor x/\mu\rfloor
}) } \\
&& \qquad \leq C \sup_{0 < \rho< \rho(x) } \rho^{(x-x^\gamma)/\mu- 1} |\log
\rho| x^{\upsilon+\alpha- 1 + \epsilon} \\
&& \qquad \leq C \sup_{t > (x^{\upsilon} \log x)^{-1}} \exp \biggl( - \biggl( \frac
{x - x^\gamma}{\mu} -1  \biggr) t + \log t + (\upsilon+\alpha-1+\epsilon
) \log x  \biggr) \\
&& \qquad = C \exp \biggl( - \frac{x^{1-\upsilon}}{\mu\log x}  \biggl(1-\frac
{1}{x^{1-\gamma}} - \frac{\mu}{x}  \biggr) - \log\log x + (\alpha
-1+\epsilon) \log x  \biggr)\\
&& \qquad  \to0 .
\end{eqnarray*}
\upqed
\end{pf*}

We now prove Theorem \ref{T.NewMain1}.

\begin{pf*}{Proof of Theorem \ref{T.NewMain1}}
Again, we only prove the statement regarding the uniformity in $0 <
\rho< 1$, since the statement for $x > 0$ follows from Lemma \ref
{L.Equivalence} and the observation that, as $\rho\nearrow1$,
\begin{eqnarray*}
&&\sup_{0< x < (1-\rho)^{-1/4}} \frac{\rho}{1-\rho}  \biggl(1-\rho^{x/\mu}
- \frac{(1-\rho) x}{\mu} \rho^{x/\mu}  \biggr) P(X_1 > x) \\
&& \qquad \leq\sup_{0< x < (1-\rho)^{-1/4}} \frac{\rho}{1-\rho}  \biggl( |\log
\rho| \frac{x}{\mu} - \frac{(1-\rho) x}{\mu} \rho^{x/\mu}  \biggr) \\
&& \qquad = \frac{\rho}{\mu(1-\rho)^{5/4}}  \bigl( |\log\rho| - (1-\rho) e^{-|\log
\rho| (1-\rho)^{-1/4}/\mu}  \bigr) \\
&& \qquad = O  \bigl( (1-\rho)^{1/2}  \bigr) .
\end{eqnarray*}

By Theorem \ref{T.SumsApprox} we only need to show that
\[
\sup_{0 < \rho< 1}  \biggl| \frac{S(\rho,x) -  ({\rho}/({1-\rho})) \gamma
(x,\rho) P(X_1 > x)}{ ({\rho}/({1-\rho})) \gamma(x,\rho) P(X_1 > x) +
\rho^{x/\mu}}  \biggr| \to0
\]
as $x \to\infty$. We start by noting that
\begin{eqnarray*}
&&S(\rho,x) - \frac{\rho}{1-\rho} \gamma(x,\rho) P(X_1 > x) \\
&& \qquad = \sum_{n=1}^{M(x)} (1-\rho) \rho^n n  \bigl( P\bigl(X_1 > x - (n-1)\mu\bigr) -
P(X_1 > x )  \bigr) \\
&& \qquad \quad  {}  + P(X_1 > x)  \Biggl( \sum_{n=1}^{M(x)} (1-\rho) \rho^n n -
\frac{\rho}{1-\rho} \gamma(x,\rho)  \Biggr) .
\end{eqnarray*}
Then, since $ \frac{\rho}{1-\rho} \gamma(x,\rho) \geq\sum
_{n=1}^{\lfloor x/\mu\rfloor-1} (1-\rho) \rho^n n$,
\begin{eqnarray*}
&& \biggl| S(\rho,x) - \frac{\rho}{1-\rho} \gamma(x,\rho) P(X_1 > x)
\biggr| \\
&& \qquad \leq P(X_1 > x) \sum_{n=1}^{M(x)} (1-\rho) \rho^n n  \biggl( \frac{P(X_1
> x - (n-1)\mu)}{P(X_1 > x)} - 1  \biggr) \\
&& \qquad \quad  {}  + P(X_1 > x) \Biggl ( \frac{\rho}{1-\rho} \gamma(x,\rho) -
\sum_{n = 1}^{M(x)} (1-\rho) \rho^n n  \Biggr) .
\end{eqnarray*}
The second term can be bounded as follows:\vspace*{-1.5pt}
%
%e3.6 ###
\begin{eqnarray} \label{eq:EasyTerm}
&&P(X_1 > x)  \Biggl( \frac{\rho}{1-\rho} \gamma(x,\rho) - \sum_{n =
1}^{M(x)} (1-\rho) \rho^n n  \Biggr) \notag\\
&& \qquad \leq\frac{\rho}{1-\rho} P(X_1 > x) \bigl(\rho^{M(x)} - \rho^{x/\mu}\bigr)\biggl (
1 + \frac{(1-\rho)x}{\mu}  \biggr) \\
&& \qquad \leq\frac{C}{1-\rho} P(X_1 > x)\bigl (\rho^{(x-x^\beta)/\mu} - \rho^{x/\mu
+1}\bigr)\bigl  ( 1 + (1-\rho)x  \bigr) ,\notag
\end{eqnarray}
where $C > 0$ is a generic constant. Fix $0 < \epsilon< \beta(\alpha
-2) /(\alpha-2+\beta)$ and set $N(x) = \lfloor(1-\epsilon) x/\mu
\rfloor$. Then, for $1 \leq n \leq N(x)$,\vspace*{-1.5pt}
\begin{eqnarray*}
\frac{P(X_1 > x - (n-1)\mu)}{P(X_1 > x)} &\leq& \biggl( \frac{x-(n-1)\mu
}{x}  \biggr)^{-\alpha+1} \sup_{1 \leq n \leq N(x)} \frac{L(x-(n-1)\mu
)}{L(x)} \\
&\leq& \biggl( 1 + (\alpha-1) \epsilon^{-\alpha-2} \frac{(n-1)\mu}{x}
 \biggr) \bigl(1+ \varphi_1(x)\bigr) ,
\end{eqnarray*}
where $\varphi_1(x) = \sup_{\epsilon\leq t \leq1} \frac{L(tx)}{L(x)}
- 1 \to0$ by properties of slowly varying functions. Therefore, for $1
\leq n \leq N(x)$,\vspace*{-1pt}
\[
\frac{P(X_1 > x - (n-1)\mu)}{P(X_1 > x)} -1 \leq C  \biggl( \frac{n-1}{x}
+ \varphi_1(x)  \biggr) .
\]
It follows that\vspace*{-1.5pt}
%
%e3.7 ###
\begin{eqnarray}\label{eq:UpToN}
&&P(X_1 > x) \sum_{n=1}^{N(x)} (1-\rho) \rho^n n  \biggl( \frac{P(X_1 > x
- (n-1)\mu)}{P(X_1 > x)} - 1  \biggr) \notag\\
&& \qquad \leq C P(X_1 > x)  \Biggl( \frac{1}{x} \sum_{n=1}^{N(x)} (1-\rho) \rho^n
n(n-1) + \varphi_1(x) \frac{\rho}{1-\rho} \gamma(x,\rho)  \Biggr) \notag
\\
&& \qquad \leq C P(X_1 > x)  \biggl( \frac{2\rho^2}{x (1-\rho)^2}  \bigl(1 - \rho
^{N(x)} - (1-\rho)N(x)\rho^{N(x)}  \bigr)   \\
&&  \qquad  \quad \hspace*{156pt}   {}   + \varphi_1(x) \frac{\rho}{1-\rho} \gamma(x,\rho)
 \biggr) \notag\\
&& \qquad \leq\frac{C\rho}{1-\rho} P(X_1 > x) \gamma(x,\rho) \biggl ( \frac{1}{x
(1-\rho)} + \varphi_1(x)  \biggr).\notag
\end{eqnarray}
For the terms $N(x) < n \leq M(x)$ we have\vspace*{-1.5pt}
%
%e3.8 ###
\begin{eqnarray}\label{eq:FromNtoM}
&&P(X_1 > x) \sum_{n=N(x)+1}^{M(x)} (1-\rho) \rho^n n \biggl ( \frac{P(X_1
> x - (n-1)\mu)}{P(X_1 > x)} - 1  \biggr) \notag\\
&& \qquad \leq C x (1-\rho) \rho^{N(x)+1} \sum_{n=N(x)+1}^{M(x)} P\bigl(X_1 >
x-(n-1)\mu\bigr) \\
&& \qquad \leq C x (1-\rho) \rho^{(1-\epsilon)x/\mu} \int_{x-\mu M(x)}^{x-\mu
N(x)} P(X_1 > t) \,dt \nonumber\\
&& \qquad \leq C x(1-\rho) \rho^{(1-\epsilon)x/\mu} \bigl(x - \mu M(x)\bigr) P\bigl(X_1 > x-\mu
M(x)\bigr) \notag\\
&& \qquad \leq C (1-\rho) \rho^{(1-\epsilon)x/\mu} x^{1+\beta} P(X_1 > x^\beta),\notag
\end{eqnarray}
where for the third inequality we used Proposition 1.5.10 in \cite{BiGoTe1987}.
Combining \eqref{eq:EasyTerm}--\eqref{eq:FromNtoM} gives
\begin{eqnarray*}
&& \biggl| S(\rho,x) - \frac{\rho}{1-\rho} (1- \rho^{x/\mu}) P(X_1 > x)
 \biggr| \\
&& \qquad \leq\frac{C}{1-\rho} P(X_1 > x)  \bigl( \rho^{(x-x^\beta)/\mu} - \rho
^{x/\mu+1} \bigr)\bigl (1 + (1-\rho)x\bigr) \\
&& \qquad  \quad {}  + \frac{C\rho}{1-\rho} P(X_1 > x) \gamma(x,\rho)  \biggl(
\frac{1}{x (1-\rho)} + \varphi_1(x)  \biggr) \\
&& \qquad  \quad {}  + C (1-\rho) \rho^{(1-\epsilon)x/\mu} x^{1+\beta} P(X_1 >
x^\beta) .
\end{eqnarray*}
Let $A(\rho,x) = \frac{\rho}{1-\rho} \gamma(x,\rho) P(X_1 > x) + \rho
^{x/\mu}$ and define $\rho(x) = 1 - c\mu(\alpha-2)\log x/x$, with $\frac
{(1-\beta)(\alpha-2+\epsilon)}{(\alpha-2) (1-\epsilon) } < c < 1$. Note
that $\gamma(x,\rho) \sim1$ as $x \to\infty$ uniformly for $0 < \rho
\leq\rho(x)$. Then,
\begin{eqnarray*}
&&\sup_{0 < \rho\leq\rho(x)} \frac{1}{A(\rho,x)}  \biggl| S(\rho,x) -
\frac{\rho}{1-\rho} (1- \rho^{x/\mu}) P(X_1 > x)  \biggr| \\
&& \qquad \leq C \sup_{0 < \rho\leq\rho(x)} \biggl \{ \frac{1}{\gamma(x,\rho)}
 \bigl( \rho^{(x-x^\beta)/\mu-1} - \rho^{x/\mu} \bigr)\bigl (1 + (1-\rho)x\bigr)
+ \frac{1}{x(1-\rho)}  \\
&&  \qquad  \quad \hspace*{96pt}{}    + \varphi_1(x) + \frac
{(1-\rho)^2 x^{1+\beta} P(X_1 > x^\beta)}{\rho\gamma(x,\rho) P(X_1 >
x)} \rho^{(1-\epsilon)x/\mu}  \biggr\} \\
&& \qquad \leq C  \biggl\{ \rho(x)^{(x-x^\beta)/\mu-1} \bigl(1 + \bigl(1-\rho(x)\bigr)x\bigr)
 + \frac{1}{x(1-\rho(x))}  + \varphi_1(x)  \\
&&  \qquad  \quad\hspace*{57pt} {}   + \frac
{P(X_1 > x^\beta)}{ P(X_1 > x)} x^{1+\beta} \bigl(1-\rho(x)\bigr)^2 \rho
(x)^{(1-\epsilon)x/\mu-1}  \biggr\} \\
&& \qquad \leq C \biggl \{ \frac{\log x}{x^{c(\alpha-2)}} + \frac{1}{\log x} +
\varphi_1(x) + \frac{P(X_1 > x^\beta)}{ P(X_1 > x)} \cdot\frac{(\log
x)^2}{x^{1-\beta+c(\alpha-2) (1-\epsilon)} }  \biggr\}.
\end{eqnarray*}
The first three terms in the expression above clearly converge to zero.
To see that the fourth one does as well use\ Potter's theorem [\cite
{BiGoTe1987}, page 25] to obtain
\begin{eqnarray*}
&&\frac{P(X_1 > x^\beta)}{ P(X_1 > x)} \cdot\frac{(\log x)^2}{x^{1-\beta
+c(\alpha-2) (1-\epsilon)} } \\
&& \qquad \leq A_\epsilon  x^{(1-\beta)(\alpha
-1+\epsilon)} \cdot\frac{(\log x)^2}{x^{1-\beta+c(\alpha-2) (1-\epsilon
)} }
\end{eqnarray*}
for some constant $A_\epsilon> 1$. Our choice of $\epsilon$ and $c$
guarantees that $1-\beta+c(\alpha-2)(1-\epsilon) > (1-\beta)(\alpha
-1+\epsilon)$.

To analyze the supremum over $\rho(x) \leq\rho< 1$ we first note that
\[
\frac{\rho}{(1-\rho)^2} \gamma(x,\rho) \leq\sum_{n=1}^{\lfloor x/\mu
\rfloor+1} \rho^n n %+ \rho^{\lfloor x/\mu\rfloor} (1-\rho)^2(\lfloor
%x/\mu\rfloor+ 1)
.
\]
Then,
\begin{eqnarray*}
&&\sup_{\rho(x) \leq\rho< 1} \frac{1}{A(\rho,x)}  \biggl| S(\rho,x) -
\frac{\rho}{1-\rho} (1- \rho^{x/\mu}) P(X_1 > x)  \biggr| \\
&& \qquad \leq C \sup_{\rho(x) \leq\rho< 1} \biggl \{\frac{P(X_1 > x)}{1-\rho}
 \bigl( \rho^{-x^\beta/\mu} - \rho \bigr) \bigl(1 + (1-\rho)x\bigr)   \\
&&  \qquad  \quad\hspace*{51pt} {}
 + \frac{\rho\gamma(x,\rho) P(X_1 > x) }{(1-\rho
) \rho^{x/\mu}}\biggl(\frac1{x(1-\rho
)}+\varphi_1(x)\biggr)
 \\
&&  \qquad  \quad\hspace*{91pt} {}  + (1-\rho) \rho^{-\epsilon x/\mu}
x^{1+\beta} P(X_1 > x^\beta)  \biggr\} \\
&& \qquad \leq C \sup_{\rho(x) \leq\rho< 1} \Biggl \{\frac{P(X_1 > x)}{1-\rho}
|\log\rho| x^\beta\log x + \rho^{-\epsilon x/\mu} x^{\beta} \log x
P(X_1 > x^\beta) \\
&&  \qquad  \quad\hspace*{101pt} {}
+ \frac{P(X_1>x)}{x}\sum^{\lceil x/\mu \rceil+1}_{n=1}
\rho^{n-x/\mu}n\bigl(1+x(1-\rho)\bigr)\Biggr\}
\\
&& \qquad \leq C  \Biggl\{ P(X_1 > x) x^\beta\log x + \rho(x)^{-\epsilon x/\mu}
x^\beta\log x P(X_1 > x^\beta)   \\
&&  \qquad  \quad \hspace*{54pt}{}
+ P(X_1>x) \sum^{\lceil x/\mu \rceil+1}_{n=1}
\rho(x)^{n-x/\mu}(1+\log x)\Biggr\} .
\end{eqnarray*}
The first term clearly converges to zero. To see that the second and third terms converge to zero as well note that
\begin{eqnarray*}
\rho(x)^{-\epsilon x/\mu} x^\beta\log x P(X_1 > x^\beta) &\leq& C
x^{\epsilon c(\alpha-2) + \beta- \beta(\alpha-1-\epsilon)} \log x \\
&\leq& C x^{\epsilon(\alpha-2) - \beta(\alpha-2-\epsilon)} \log x
\end{eqnarray*}
and
\begin{eqnarray*}
P(X_1 > x) \sum_{n=1}^{\lceil x/\mu\rceil+1} \rho(x)^{n-x/\mu}\log x &\leq& C x^{-\alpha+\epsilon} \frac{ ( 1 - \rho
(x)^{\lceil x/\mu\rceil+ 1}  )}{\rho(x)^{ x/\mu
}(1-\rho(x))} \log x\\
&\leq& C x^{-\alpha+\epsilon+1 + c(\alpha-2)}  \\
&\leq& C x^{-1+\epsilon} .
\end{eqnarray*}
Our choice of $\epsilon$ guarantees that both expressions above
converge to zero. This completes the proof.
\end{pf*}

We end this section with the proof of Corollary \ref{C.Threshold2}.

\begin{pf*}{Proof of Corollary \ref{C.Threshold2}}
Let
\[
y = y(\rho) = c \mu(\alpha-2) (1-\rho)^{-1} \log \bigl( (1-\rho)^{-1}
 \bigr) .
\]
We start with the proof of part (a). We need to verify that for
$0 < c \leq1$
\[
\sup_{0 \leq x \leq y}  \biggl| \frac{ ({\rho}/({1-\rho}))\gamma(x,\rho)
P(X_1 > x) + \rho^{x/\mu}}{\exp(-(1-\rho)x/\mu)} - 1  \biggr| \to0
\]
as $\rho\nearrow1$. Note that
%
%e3.10 ###
%e3.9 ###
\begin{eqnarray}
&&\sup_{0 \leq x \leq y}  \biggl| \frac{ ({\rho}/({1-\rho}))\gamma(x,\rho)
P(X_1 > x) + \rho^{x/\mu}}{\exp(-(1-\rho)x/\mu)} - 1  \biggr| \notag\\\label{eq:NonNegligible}
&& \qquad \leq \sup_{0 \leq x \leq y}  \biggl| \exp \biggl( \frac{x \log\rho}{\mu}
+ \frac{(1-\rho)x}{\mu}  \biggr) - 1  \biggr|  \\\label{eq:Negligible}
&&  \qquad  \quad {} + \sup_{0 \leq x \leq y} \frac{\rho\gamma(x,\rho
)}{(1-\rho)} P(X_1>x) \exp\bigl((1-\rho)x/\mu\bigr)  .
\end{eqnarray}
We can bound \eqref{eq:NonNegligible} as follows:
\begin{eqnarray*}
\sup_{0 \leq x \leq y}  \biggl| \exp \biggl( \frac{x \log\rho}{\mu} + \frac
{(1-\rho)x}{\mu}  \biggr) - 1  \biggr| &\leq&\sup_{0 \leq x \leq y} \frac
{x |\log\rho+ 1-\rho|}{\mu} \\
&\leq& C y (1-\rho)^2 \\
&\leq& C t^{-1} \log t ,
\end{eqnarray*}
where $t = (1-\rho)^{-1}$. Also, note that since $\rho^{x/\mu} \geq1 -
|\log\rho|x/\mu$ and $|\log\rho| = 1-\rho+ O((1-\rho)^2)$ as $\rho
\nearrow1$,
\begin{eqnarray*}
\gamma(x,\rho) &=& 1 - \rho^{x/\mu} - \rho^{x/\mu} (1-\rho)x/\mu\\
&\leq&\bigl(|\log\rho| - (1-\rho)\bigr) x/\mu+ |\log\rho| (1-\rho) (x/\mu)^2 \\
&\leq& C (1-\rho)^2 x^2.
\end{eqnarray*}
Then \eqref{eq:Negligible} is bounded by
\begin{eqnarray*}
&&\sup_{0 \leq x \leq y} \frac{\rho\gamma(x,\rho)}{(1-\rho)} P(X_1>x)
\exp\bigl((1-\rho)x/\mu\bigr) \\
&& \qquad \leq C \sup_{0 \leq x \leq(1-\rho)^{-1/4}} \frac{ \gamma(x,\rho)
}{1-\rho}\\
&&  \qquad  \quad {} + \sup_{(1-\rho)^{-1/4} \leq x \leq y} \frac{1}{1-\rho} P(X_1
> x) \exp \bigl((1-\rho)x/\mu \bigr) \\
&& \qquad \leq\sup_{0 \leq x \leq(1-\rho)^{-1/4}} C (1-\rho) x^2 \\
&&  \qquad  \quad {} + \sup_{(1-\rho)^{-1/4} \leq x \leq y} \frac{L(x)}{1-\rho
} \exp \bigl( (1-\rho)x/\mu- (\alpha-1)\log x  \bigr) \\
&& \qquad \leq C (1-\rho)^{1/2} + \frac{L(y)}{1-\rho} \exp \bigl( (1-\rho)y/\mu-
(\alpha-1)\log y  \bigr) \\
&& \qquad \leq C t^{-1/2} + C L(t\log t) \exp \bigl( -(1-c)(\alpha-2) \log t -
(\alpha-1) \log\log t  \bigr) .
\end{eqnarray*}
Clearly, if $0 < c < 1$, then the two expressions above converge to
zero as $t \to\infty$. If $c = 1$ and $L(x)/(\log x)^{\alpha-1} \to0$
as $x \to\infty$, then
\begin{eqnarray*}
&&\lim_{\rho\nearrow1} \sup_{0 \leq x \leq y}  \biggl| \frac{ ({\rho
}/({1-\rho})) \gamma(x,\rho) P(X_1 > x) + \rho^{x/\mu}}{\exp(-(1-\rho)x/\mu
)} - 1  \biggr| \\
&& \qquad \leq C \lim_{t\to\infty} L(t \log t) \exp \bigl( - (\alpha-1) \log\log
t  \bigr) \\
&& \qquad = C \lim_{t\to\infty} \frac{L(t \log t)}{(\log(t \log t))^{\alpha-1}}
\cdot \biggl(\frac{\log(t \log t)}{\log t}  \biggr)^{\alpha-1} = 0 .
\end{eqnarray*}
We now move to part (b). We need to verify that for $c \geq1$
\[
\sup_{x \geq y}  \biggl| \frac{( {\rho}/({1-\rho}))\gamma(x,\rho) P(X_1 >
x) + \rho^{x/\mu}}{ ({\rho}/({1-\rho})) P(X_1 > x)} - 1  \biggr| \to0
\]
as $\rho\nearrow1$. Note that
\begin{eqnarray*}
&&\sup_{x \geq y}  \biggl| \frac{ ({\rho}/({1-\rho})) \gamma(x,\rho) P(X_1
> x) + \rho^{x/\mu}}{ ({\rho}/({1-\rho})) P(X_1 > x)} - 1  \biggr| \\
&& \qquad \leq\sup_{x \geq y}  | \gamma(x,\rho) - 1  | + \sup_{x \geq
y} \frac{(1-\rho) \rho^{x/\mu}}{\rho P(X_1 > x)} \\
&& \qquad \leq\rho^{y/\mu} \bigl(1 + (1-\rho)y/\mu\bigr) + C \sup_{x \geq y} \frac{1-\rho
}{L(x)} \exp \biggl( -\frac{x}{\mu} (1-\rho) + (\alpha-1)\log x  \biggr) \\
&& \qquad \leq C \rho^{y/\mu}(1-\rho) y + C \frac{1-\rho}{L(y)} \exp \biggl( -\frac
{y}{\mu} (1-\rho) + (\alpha-1)\log y  \biggr) \\
&& \qquad \leq C t^{- c(\alpha-2)} \log t\\
&& \qquad  \quad {} + \frac{C}{L(t\log t)} \exp \bigl( -
(c-1) (\alpha-2) \log t + (\alpha-1) \log\log t  \bigr) ,
\end{eqnarray*}
where $t = (1-\rho)^{-1}$. Clearly, if $c > 1$ the above converges to
zero. If $c = 1$ and $L(x)/(\log x)^{\alpha-1} \to\infty$ as $x \to
\infty$, then
\begin{eqnarray*}
&&\lim_{\rho\nearrow1} \sup_{0 \leq x \leq y}  \biggl| \frac{ ({\rho
}/({1-\rho}))\gamma(x,\rho) P(X_1 > x) + \rho^{x/\mu}}{\exp(-(1-\rho)x/\mu
)} - 1  \biggr| \\
&& \qquad \leq C \lim_{t \to\infty} \frac{1}{L(t\log t)} \exp \bigl( (\alpha-1)
\log\log t  \bigr) \\
&& \qquad = C \lim_{t\to\infty} \frac{(\log(t \log t))^{\alpha-1}}{L(t \log t)}
\cdot \biggl(\frac{\log t}{\log(t \log t)}  \biggr)^{\alpha-1} = 0 .
\end{eqnarray*}
\upqed
\end{pf*}

%%% Numerical Approximations
%s4 ###
\section{Numerical approximations} \label{S.Numerical}

Theorems \ref{T.SumsApprox} and \ref{T.NewMain1} suggest approximating
$P(W_\infty(\rho) > x)$ either with
\[
H(\rho,x) \triangleq S(\rho,x) + \rho^{x/\mu}
= \sum_{n=1}^{M(x)} (1-\rho)\rho^n n P\bigl(X_1 > x-(n-1)\mu\bigr) + \rho^{x/\mu}
\]
or with
\[
J(\rho,x) \triangleq\frac{\rho}{1-\rho} \gamma(\rho,x) P(X_1 > x) +
\rho^{x/\mu},
\]
respectively.

We compared both approximations to simulated values of $P(W_\infty(\rho
) > x)$ and found that $H(\rho,x)$ tends to be better than $J(\rho,x)$
and seems to perform very well across all values of $x$ for different
choices of $\rho$. This is not surprising given that $H(\rho,x)$ more
closely resembles the   Pollaczek--Khintchine formula than $J(\rho,x)$.

When $\sigma^2 = \var(X_1) < \infty$, the central limit theorem can be
used to approximate the tail of the Pollaczek--Khintchine formula in a
way that $\sigma^2$ is incorporated into the approximation. The term
$\rho^{x/\mu}$ appearing in the definitions of $H(\rho,x)$ and $J(\rho
,x)$ can be replaced by
\[
T(\rho, x) \triangleq\sum_{n=1}^\infty(1-\rho) \rho^n  \bigl( 1- \Phi
 \bigl((x-n\mu)/\sqrt{\sigma^2 n} \bigr)  \bigr) ,
\]
which can alternatively be written as $T(\rho,x) = E [ \rho
^{M(x,Z)}  ]$, where $Z \sim$ $N(0,1)$ and $M(x,z) = \lfloor (
\sqrt{x/\mu+(\sigma z)^2/(2\mu)^2} - (\sigma z)/(2\mu)  )^2
\rfloor$. We do not give proofs here, but it can be shown that provided
$\sigma^2 < \infty$, Theorems \ref{T.SumsApprox} and \ref{T.NewMain1}
continue to hold with $\rho^{x/\mu}$ replaced by $T(\rho,x)$. This is
relevant from the numerical standpoint since the resulting
approximations tend to perform better than those with the simpler $\rho
^{x/\mu}$.

%f1 ###
\begin{figure}

\includegraphics{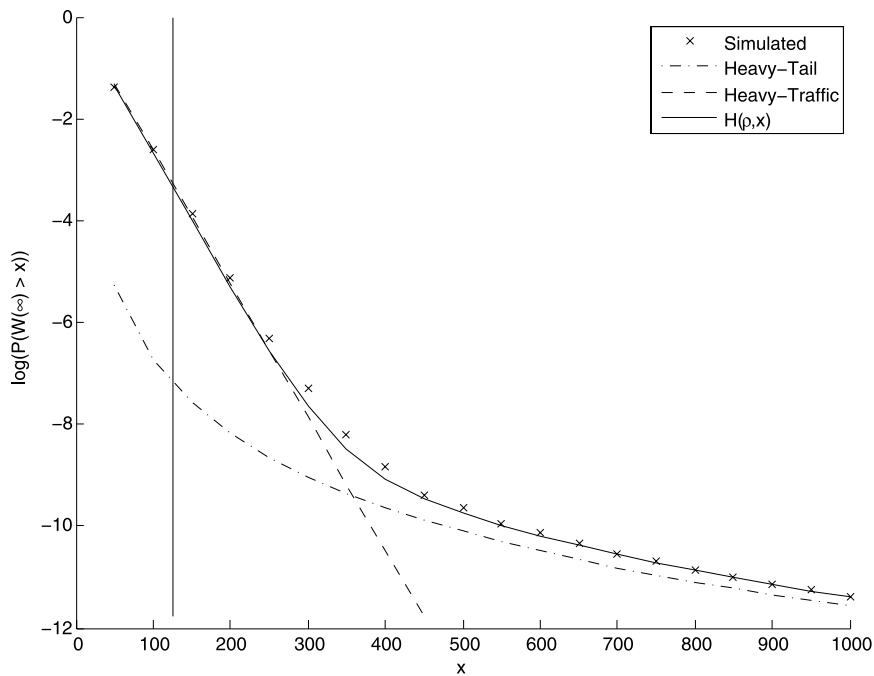}

\caption{Pareto integrated tail with $\rho= 0.95$ and $\alpha= 3.1$.}\label{fig1}
\end{figure}

%f2 ###
\begin{figure}

\includegraphics{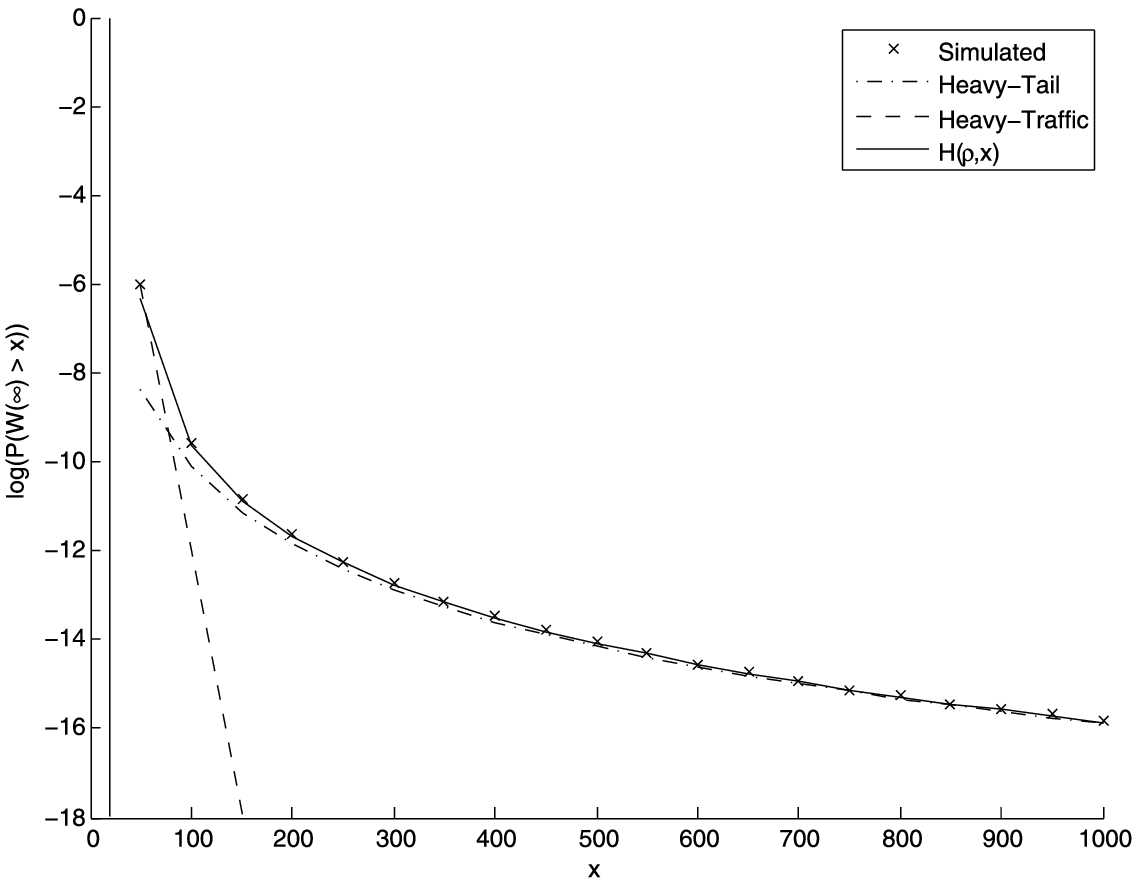}

\caption{Pareto integrated tail with $\rho= 0.8$ and $\alpha= 3.5$.}\label{fig2}
\end{figure}

We plotted approximation $H(\rho,x)$ against simulated values of
$P(W_\infty(\rho) > x)$. Figures~\ref{fig1} and~\ref{fig2} correspond to queues having
Pareto integrated tail distribution, that is, $P(X_1 > x) = x^{-\alpha
+1}$ for $x \geq1$. For comparison purposes we also plotted the
heavy-traffic approximation,
\[
\mathit{Heavy\mbox{-}Traffic} = \exp\bigl(-(1-\rho)x/\mu\bigr),
\]
and the heavy-tail asymptotic,
\[
\mathit{Heavy\mbox{-}Tail} = \frac{\rho}{1-\rho} P(X_1 > x) .
\]
The vertical line corresponds to the value
\[
\hat x(\rho) = \mu(\alpha-2) (1-\rho)^{-1} \log \bigl((1-\rho)^{-1}
\bigr) .
\]

The simulated values of $P(W_\infty(\rho) > x)$ were obtained using the
conditional Monte Carlo algorithm from \citet{AsmKro06}, and each point
was estimated using enough simulation runs to obtain a relative error
of at most 0.05 with approximately $99\%$ confidence.

% imsref loaded by smiklovaite, 2010-09-15 15:38:40
%

\printaddresses

\end{document}